\documentclass[10pt]{article}
\usepackage{color}
\usepackage{tikz}
\usepackage{amssymb, amsmath}
\usepackage{epsfig}
\usepackage{graphicx}

\newcommand{\bbc}{{\mathcal C}}

 \newcommand{\bbR}{{\mathbb R}}

\newtheorem{theorem}{Theorem}
 
 \newtheorem{proposition}{Proposition}

\newcommand{\proof}{\noindent\textbf{Proof.~}}
\newcommand{\vep}{\varepsilon}
\newcommand{\qed}{\space\hfill\hspace*{\fill} $\vbox{\hrule\hbox{\vrule
height1.3ex\hskip1.3ex\vrule}\hrule}$\hss\vskip\topsep\relax}

\begin{document}

\title{On construction of boundary preserving numerical schemes}

\author{Nikolaos Halidias \\
{\small\textsl{Department of Mathematics }}\\
{\small\textsl{University of the Aegean }}\\
{\small\textsl{Karlovassi  83200  Samos, Greece} }\\
{\small\textsl{email: nikoshalidias@hotmail.com}}}

\maketitle

\begin{abstract}Our aim in this note is to extend the semi discrete technique by combine it with the split step method.  We apply
our new method to the Ait-Sahalia model and propose an explicit
and positivity preserving numerical scheme.
\end{abstract}



\section{Introduction}
In this paper we describe a technique to construct numerical
schemes by combining the split step method (see for example
\cite{schurz1}) and the semi discrete method (see
\cite{Halidias1}, \cite{Halidias2}, \cite{Halidias3}). Using the
semi discrete method we have constructed explicit and positivity
numerical schemes for various stochastic differential equations
arising in finance (see \cite{Halidias4}, \cite{Halidias5},
\cite{Halidias6}, \cite{Halidias7}, \cite{Halidias8}).

Using the proposed technique (split-step and semi-discrete) we are
able to handle more situations in which we want to construct
explicit and boundary preserving numerical schemes. Our starting
point was the paper \cite{Mao} (see also \cite{Andreas}) in which
the authors proposed an implicit numerical scheme to approximate
the Ait-Sahalia model (see \cite{Ait}), which is the following,
\begin{eqnarray*}
x_t = x_0 + \int_0^t \left(\frac{a_1}{x_s} - a_2 + a_3x_s -
a_4x_s^{r} \right) ds + \int_0^t \sigma x_s^{\rho} dw_s
\end{eqnarray*}
with $x_0 \in \mathbb{R}_+$. We assume that all the coefficients
are nonnegative and that $r+1 > 2 \rho$.

Using this new method we describe and analyze a new explicit and
positivity preserving numerical scheme for the Ait-Sahalia model
which arise in finance. As far as we know this is the first
explicit scheme for this model, however this does not mean that
from the computational point of view is cheaper than the implicit
ones (\cite{Mao}, \cite{Andreas}). We have to   make extended
numerical experiments in order to compare them.

\section{The drift splitting}
Let $(\Omega, {\cal F}, \mathbb{P}, {\cal F}_t)$ be a complete
probability space with a filtration and let  a Wiener process
$(W_t)_{t \geq 0}$ defined on this space. Consider the following
stochastic differential equation,
\begin{eqnarray}
x_t = x_0 + \int_0^t a(s,x_s)ds + \int_0^t b(s,x_s)dW_s,
\end{eqnarray}
where $a,b: \mathbb{R}_+ \times \mathbb{R} \to \mathbb{R}$ are
measurable functions
 and $x_0$ such that is ${\cal F}_0$-measurable and square
 integrable.

Let $0 = t_0 < t_1 < ...<t_n = T$ and set $\Delta = \frac{T}{n}$.
In general, one can split the sde above in $m$ equations. For
example, if $a(t,x) = a_1(t,x) + ... +a_m(t,x)$  then we can have
the following splitting
\begin{eqnarray*}
y_m(0) & = & x_0 \\
y_1(t) &  = & y_m(t_n) + \int_{t_n}^t a_1(s,y_1(s))ds, \quad t \in (t_n,t_{n+1}] \\
& \vdots  & \\ y_{m-1}(t) & = &  y_{m-2}(t_{n+1})+\int_{t_n}^t
a_{m-1}(s,y_{m-1}(s))ds, \quad t \in (t_n,t_{n+1}] \\  y_m(t) & =
& y_{m-1}(t_{n+1}) + \int_{t_n}^t a_m(s,y_m(s))ds + \int_{t_n}^t
b(s,y_m(s)) dw_s, \quad t \in (t_n,t_{n+1}]
\end{eqnarray*}
and approximate each of the above equations by a semi discrete
scheme (or another converging scheme).

Then, we can write, for $t \in (t_n,t_{n+1}]$ and $y_m(0) = x_0$
\begin{eqnarray*}
y_m(t) = y_m(t_n) & +  & \int_{t_n}^{t_{n+1}}
\left(a_1(s,y_1(s))+...+a_{m-1}(s,y_{m-1}(s))\right) ds \\ & + &
\int_{t_n}^t a_m(s,y_m(s)) ds   +   \int_{t_n}^t b(s,y_m(s)) dw_s
\end{eqnarray*}
or
\begin{eqnarray*}
y_m(t) = x_0 & +  & \int_{0}^{t_{n+1}}
\left(a_1(s,y_1(s))+...+a_{m-1}(s,y_{m-1}(s))\right) ds \\ & + &
\int_{0}^t a_m(s,y_m(s)) ds  + \int_{0}^t b(s,y_m(s)) dw_s
\end{eqnarray*}

We shall denote by $\hat{t} = t_n$ when $ t \in [t_n,t_{n+1}]$ and
$\tilde{t} = t_{n+1}$ when $t \in [t_n,t_{n+1}]$.

\section{Approximating with the semi discrete method}
Suppose that there are functions
$f_1(t,x,y,z)$,...,$f_m(t,x,y,z)$,$f_{m+1}(t,x,y)$
 such that $f_i(t,x,x,x) = a_i(t,x)$ for $i=1,...,m$ and $f_{m+1}(t,x,x) =
 b(t,x)$.

Our numerical scheme depends on the choices of $f_i$ and therefore
we should impose conditions on them. For fixed $a(t,x),b(t,x)$ one
can choose different $f_i$ in such a way that the corresponding
numerical schemes does not converge.

Denoting again our approximation by $y_m(t)$ we write, for $t \in
(t_n,t_{n+1}]$
\begin{eqnarray*}
y_m(0) & = & x_0, \nonumber \\
y_1(t) & = & y_m(t_n) + \int_{t_n}^t
f_1(s,y_m(t_n)),y_1(t),y_1(s))ds,
 \\
y_2(t) & = & y_1(t_{n+1}) + \int_{t_n}^t
f_2(s,y_m(t_n),y_2(t),y_2(s))ds,  \\
& \vdots & \\
y_{m-1}(t) & = & y_{m-2}(t_{n+1}) + \int_{t_n}^t
f_{m-1}(s,y_m(t_n),y_{m-1}(t),y_{m-1}(s))ds,  \\
y_m(t) & = & y_{m-1}(t_{n+1}) +  \int_{t_n}^t
f_{m}(s,y_m(t_n),y_{m}(s))ds + \int_{t_n}^t b(s,y_m(t_n),y_m(s))
dw_s,
\end{eqnarray*}
It is obvious that we should choose $f_i$ in such a way that all
the above equations has at least one strong solution. Then, we
have constructed an approximation scheme which is $y_m(t)$ and
under suitable conditions we will show that this converges
strongly to the unique strong solution of problem (1). If some of
the above equations admits more than one solution then we have
constructed at least two approximation schemes and we choose the
numerical scheme that has  the desired properties, positivity
preserving for example.

In a more compact form we can write, for $t \in (t_n,t_{n+1}]$,

\begin{eqnarray}
y_m(t) = y_m(t_n) & +  & \int_{t_n}^{t_{n+1}}
\Big(f_1(s,y_1(s),y_1(t),y_m(t_n))+...+f_{m-1}(s,y_{m-1}(s),y_{m-1}(t),y_m(t_n))\Big)
ds \nonumber \\ & + & \int_{t_n}^t f_m(s,y_m(s), y_m(t_n)) ds
  +  \int_{t_n}^t f_{m+1}(s,y_m(s),y_m(t_n)) dw_s, \; t
\in (t_n,t_{n+1}],
\end{eqnarray}
with $y_m(0) = x_0$, and also
\begin{eqnarray}
y_m(t) = x_0 & +  & \int_{0}^{t_{n+1}}
\Big(f_1(s,y_1(s),y_1(t),y_m(\hat{s}))+...+f_{m-1}(s,y_{m}(s),y_{m-1}(t),y_{m}(\hat{s}))\Big)
ds \nonumber \\ & + & \int_{0}^t f_m(s,y_m(s), y_m(\hat{s})) ds  +
 \int_{0}^t f_{m+1}(s,y_m(s),y_m(\hat{s})) dw_s, \quad  t \in
(t_n,t_{n+1}]
\end{eqnarray}

Furthermore we have, for $i=1,...,m-1$ and for $t \in
(t_n,t_{n+1}]$
\begin{eqnarray}
y_i(t)  =  y_i(t_n) & + & \int_{t_n}^{t_{n+1}}
\Big(f_1(s,y_1(s),y_1(t),y_m(\hat{s}))+...+f_{i-1}(s,y_{i-1}(s),y_{i-1}(t),y_{m}(\hat{s}))\Big)
ds \nonumber \\ &  + & \int_{t_{n-1}}^{t_n}
\Big(f_{i+1}(s,y_{i+1}(s),y_{i+1}(t),y_{m}(\hat{s}))+...+f_{m}(s,y_{m}(s),y_{m}(t),y_{m}(\hat{s}))\Big)
ds \nonumber \\ &  +  & \int_{t_n}^t f_i
(s,y_{i}(s),y_{i}(t),y_{m}(\hat{s}))ds + \int_{t_{n-1}}^{t_n}
f_{m+1}(s,y_{m}(s),y_{m}(\hat{s})) dw_s,
\\
y_m(t)  =  y_i(t_{n+1})  & + & \int_{t_n}^{t_{n+1}}
f_i(s,y_i(s),y_{i}(t),y_m(\hat{s}))+...+f_{m-1}(s,y_{m-1}(s),y_{m-1}(t),y_{m}(\hat{s}))
ds \nonumber \\ & + & \int_{t_n}^t
f_m(s,y_m(s),y_{m}(t),y_m(\hat{s})) ds + \int_{t_n}^t
f_{m+1}(s,y_{m}(s),y_{m}(\hat{s}))dw_s
\end{eqnarray}

 {\bf Assumption A} Suppose that problem (1)
has a unique strong solution and  that the following moment bounds
holds, for every $p > 0$,
\begin{eqnarray*}
\mathbb{E}( |y_1(t)|^p + |y_2(t)|^p +...+y_m(t)+ |x(t)|^p ) <
\infty
\end{eqnarray*}

{\bf Assumption B} Suppose that the functions $f_i$ for
$i=1,...,m$ satisfy the following locally Lipschitz condition,
\begin{eqnarray*}
 |f_i(t,x,x,x) - f_i(t,x_1,x_2,x_3)|  & \leq &  C_R (|x-x_1| +
|x-x_2|+|x-x_3| ), \quad i=2,...,m \\
|f_{m+1}(t,x,x) - f_{m+1}(t,x_1,x_2)| & \leq & C_R (|x-x_1| +
|x-x_2| + |x-x_1|^a)
\end{eqnarray*}
for $|x_1|, |x_2|, |x_3|, |x| \leq R$,  and for some $a \in
[\frac{1}{2},1)$.

\begin{theorem} Under Assumptions  A and B we have that
\begin{eqnarray*}
 \mathbb{E}  |y_m(t) - x(t)|^2  \leq C_R \Delta +
\frac{2^{p+1} \delta A}{p} + \frac{(p-2)2A}{p
\delta^{\frac{2}{p-2}}R^p} + e_{q-1} + \frac{C_R\Delta}{q e_q}
\end{eqnarray*}
where $e_q=e^{-q(q+1)/2}$ for every $q \in \mathbb{N}$. Therefore,
for every $\vep > 0$, we can  fix first big enough $q$, then small
enough $\delta$ and big enough $R$ and finally for small enough
$\Delta$ we obtain that $\mathbb{E}  |y_m(t) - x(t)|^2  \leq
\vep$.
\end{theorem}

\proof Set $ \rho_R = \inf \{ t \in [0,T]:
 |x(t)| \geq R \}$, $ \tau^i_R = \inf \{ t \in [0,T] : |y_i(t)| \geq R
\}$ for $i=1,...,m$. Let $\theta_R = \min \{ \tau^i_R,  \rho_R
\}$.

We can prove that $\mathbb{P}(\tau^{i}_r \leq T \mbox{ or } \rho_R
\leq T ) \leq \frac{2 A}{R^p}$. Using  Young inequality we obtain,
for any $\delta > 0$,
\begin{eqnarray*}
\mathbb{E}\left(\sup_{0 \leq t \leq T } |y_m(t) -x(t)|^2\right)
\leq \mathbb{E}\left( \sup_{ 0 \leq t \leq T} \left|y_m(t
\wedge\theta_R) - x(t \wedge\theta_R)\right|^2\right) +
\frac{2^{p+1} \delta A}{p} + \frac{(p-2)2A}{p
\delta^{\frac{2}{p-2}}R^p}.
\end{eqnarray*}

The difference $x(t) - y_m(t)$ is as follows,
\begin{eqnarray*}
x(t) - y_m(t) = & &  \int_0^t \sum_{i=1}^m \big(
f_i(s,x(s),x(s),x(s)) - f_i(s,y_i(s),y_i(t),y_m(\hat{s}) )\big)ds
\\ & + & \int_0^t \big(f_{m+1} \big(s,x(s),x(s)) -
f_{m+1}(s,y_m(s),y_m(\hat{s})) \big) dw_s \\ & + &
\int_t^{t_{n+1}} \sum_{i=1}^{m-1} \big( f_i(s,x(s),x(s),x(s)) -
f_i(s,y_i(s),y_i(t),y_m(\hat{s}) \big)ds
\end{eqnarray*}

We shall estimate  the term $|x(t\wedge \theta_R) - y_m(t\wedge
\theta_R)|^2$ as follows, using  Assumption  B
\begin{eqnarray}
  & & \mathbb{E} |x(t\wedge \theta_R) - y_m(t\wedge \theta_R)|^2
  \nonumber
  \\
  & \leq & C_R\Delta^2+C_R \int_0^{t \wedge \theta_R}    \sum_{i=1}^{m}
\mathbb{E}\big(|x(s) - y_m(s)|^2 +|y_m(s)-y_i(\tilde{s})|^2
\nonumber \\ & &  + |y_m(s)-y_m(\hat{s})|^2 +
|y_i(\tilde{s})-y_i(s)|^2  + |x(s)-y_m(s)|^{2a}\big) \nonumber \\
& \leq & C_R\sqrt{\Delta} + C_R \int_0^{t \wedge \theta_R}
(\mathbb{E}|x(s) - y_m(s)|^2 +\mathbb{E}|x(s) - y_m(s)|) ds
\end{eqnarray}
The term $C_R \Delta^2$ comes from the estimation of
\begin{eqnarray*}
\int_t^{t_{n+1}} \sum_{i=1}^{m-1} \big( f_i(s,x(s),x(s),x(s)) -
f_i(s,y_i(s),y_i(t),y_m(\hat{s}) \big)ds
\end{eqnarray*}

To get (6) above we have used the following, for $i=1,...,m+1$
\begin{eqnarray*}
|x(s) - y_i(s)| & \leq &  |x(s) - y_m(s)| + |y_m(s)-y_i(s)| \\ &
\leq & |x(s) - y_m(s)| + |y_m(s)-y_{i}(\tilde{s})| + |
y_{i}(\tilde{s}) - y_i(s)|
\end{eqnarray*}
combined with (4) and (5). Furthermore we have used that  for $2a
\geq 1$ it holds that
\begin{eqnarray*}
\mathbb{E}|x(s)-y_i(s)|^{2a} = \mathbb{E}|x(s)-y_i(s)|
|x(s)-y_i(s)|^{2a-1} \leq C_R \mathbb{E}|x(s)-y_i(s)|
\end{eqnarray*}

From (4) and (5) we have, for $s \in (t_n , t_{n+1} ]$ and $
i=1,...,m$,
\begin{eqnarray*}
\mathbb{E}|y_i(s \wedge \theta_R) - y_i(t_n \wedge
\theta_R)|^2 & \leq & C_R \Delta, \quad \\
\mathbb{E}|y_m(s \wedge \theta_R) - y_i(t_{n+1} \wedge
\theta_R)|^2 & \leq & C_R \Delta, \\
\mathbb{E} |y_i(t_{n+1} \wedge \theta_R) - y_i(s \wedge
\theta_R)|^2 & \leq & C_R \Delta
\end{eqnarray*}
and all these estimates are useful to get (6).

We should estimate the following quantity (and substitute this
estimation to (6)),
\begin{eqnarray*}
\mathbb{E} |x(s)-y_m(s)| \\
\end{eqnarray*}

Let the non increasing sequence $\{e_q\}_{q\in\mathbb{N}}$ with
$e_q=e^{-q(q+1)/2}$ and $e_0=1.$ We introduce the following
sequence of smooth approximations of $|x|,$
$$
\phi_q(x)=\int_0^{|x|}dy\int_0^{y}\psi_q(u)du,
$$
where the existence of the continuous function $\psi_q(u)$ with
$0\leq \psi_q(u) \leq 2/(qu)$ and support in $(e_q,e_{q-1})$ is
justified by $\int_{e_q}^{e_{q-1}}(du/u)=q.$ The following
relations hold for $\phi_q\in\bbc^2(\bbR,\bbR)$ with
$\phi_q(0)=0,$
 $$
 |x| - e_{q-1}\leq\phi_q(x)\leq |x|, \quad |\phi_{q}^{\prime}(x)|\leq1, \quad x\in\bbR, $$
 $$
 |\phi_{q}^{\prime \prime }(x)|\leq\frac{2}{q|x|}, \,\hbox{ when }  \,e_q<|x|<e_{q-1} \,\hbox{ and }  \,  |\phi_{q}^{\prime \prime }(x)|=0 \,\hbox{ otherwise. }
 $$

Applying Ito's formula on $\phi_q(x(t)-y_m(t))$ for $t \in [0,t
\wedge \theta_R]$ we get
\begin{eqnarray*}
\mathbb{E}\phi_q(x(t) - y_m(t)| & \leq & C_R \Delta +
\frac{C_R\Delta}{q e_q}+C_R\int_0^t \mathbb{E}\phi_q^{'}(x(t) -
y_m(t)) |x(s)-y_m(s)|ds \\  & &  + C_R\int_0^t
\mathbb{E}|x(s)-y_m(s)|ds
\end{eqnarray*}
Therefore
\begin{eqnarray*}
\mathbb{E}|x(t)-y_m(t)| \leq e_{q-1} + C_R \Delta+
\frac{C_R\Delta}{q e_q} + C_R\int_0^t \mathbb{E}|x(s)-y_m(s)|ds
\end{eqnarray*}
Applying Gronwall's inequality and substituting in (6) and then
again Gronwall inequality we get the desired result. \qed

An example, is the following stochastic differential equation,
\begin{eqnarray*}
x_t = x_0 + \int_0^t k(l-x_s) - dx_s^2 ds + \sigma \int_0^t
\sqrt{x_s} dw_s
\end{eqnarray*}
For this sde, we propose the following splitting, for $t \in
(t_n,t_{n+1}]$
\begin{eqnarray*}
y_2(0) & = & x_0, \\
y_1(t) & = & y_2(t_n) + \int_{t_n}^t -ky_1(s) - dy_1^2(s) ds, \\
y_2(t) & = & y_1(t_{n+1}) + \int_{t_n}^t kl ds + \sigma
\int_{t_n}^t \sqrt{y_2(s)}dw_s
\end{eqnarray*}

The first equation can be approximated as follows, denoting again
by $y_1(t)$ the approximation
\begin{eqnarray*}
y_1(t) = y_2(t_n) + \int_{t_n}^t -k y_1(s) - dy_1(s) y_2(\hat{s})
ds
\end{eqnarray*}
which produces a positive solution whenever $y_2(t_n) > 0$. The
second equation can be approximated in the spirit of
\cite{Halidias2}.

\section{Application to the Ait-Sahalia model}
In the previous section we have described a new technique to
construct numerical schemes by combining the splitting technique
and the semi discrete method. We have proved a convergence result
when the numerical scheme satisfy some classic hypotheses. Below
we shall use this technique to construct an explicit and
positivity preserving numerical scheme for the Ait-Sahalia model
which is the following,

\begin{eqnarray*}
x_t = x_0 + \int_0^t \frac{1}{2}\left(\frac{a_1}{x_s} - a_2 +
a_3x_s - a_4x_s^{r} \right) ds + \int_0^t \sigma x_s^{\rho} dw_s
\end{eqnarray*}
with $x_0 \in \mathbb{R}_+$. We first use the transformation
\begin{eqnarray*}
y_t = x_t^2
\end{eqnarray*}
By using Ito's formula we obtain
\begin{eqnarray*}
y(t) = y(0) + \int_0^t \Big(a_1 - a_2 \sqrt{y(s)} + a_3 y(s) - a_4
y^{\frac{r+1}{2}}(s) + \sigma^2 y^{\rho}(s) \Big)ds + \int_0^t 2
\sigma y^{\frac{\rho+1}{2}}(s) dw_s
\end{eqnarray*}
 We assume that  that $r+1
> 2 \rho$. Then, we split as follows, introducing a free parameter
$a
> 0$,
\begin{eqnarray}
y_2(0) & = & x_0^2 \nonumber \\
y_1(t) & = & y_2(t_n) + \int_{t_n}^t \big( a y_1(s) - a_2
\sqrt{y_1(s)} \big)ds, \\
y_2(t) & = & y_1(t_{n+1}) + \int_{t_n}^t \big( a_1 +(a_3-a)y_2(s)-
a_4y^{\frac{r+1}{2}}_2(s) + \sigma^2 y^{\rho}_2(s) \big)ds +
\int_{t_n}^t 2 \sigma y^{\frac{\rho+1}{2}}_2(s) dw_s
\end{eqnarray}

It is easy to verify that equation (7) has at least one solution
in each interval and one of  them is the following
\begin{eqnarray*}
y_1(t) = \Big(\frac{a_2}{a} + (\sqrt{y_2(t_n)} - \frac{a_2}{a})
e^{a\frac{(t-t_n)}{2}} \Big)^2
\end{eqnarray*}
which is always positive and well posed whenever $y_2(t_n) \geq
0$.

We can approximate equation (8) by using a semi discrete approach,
namely
\begin{eqnarray*}
y_2(t) = y_1(t_{n+1}) + \int_{t_n}^t a_1+y_2(s)\big(a_3-a+\sigma^2
y_2^{\rho-1}(\hat{s})-a_4y_2^{\frac{r-1}{2}}(\hat{s}) \big) ds +
\int_{t_n}^t 2 \sigma y_2(s) y_2^{\frac{\rho-1}{2}}(\hat{s}) dw_s
\end{eqnarray*}
which have a positive and known strong solution whenever
$y_1(t_{n+1}) \geq 0$. We denote again by $y_2(t)$ the
approximation of (8).

We will use later on the following forms of $y_1,y_2$, for $t \in
(t_n,t_{n+1}]$,
\begin{eqnarray*}
y_1(t)  =  y_1(t_n) & + & \int_{t_{n-1}}^{t_n}
a_1+y_2(s)\big(a_3-a+\sigma^2
y_2^{\rho-1}(\hat{s})-a_4y_2^{\frac{r-1}{2}}(\hat{s}) \big) ds \\
&   + & \int_{t_{n-1}}^{t_n} 2 \sigma y_2(s)
y_2^{\frac{\rho-1}{2}}(\hat{s}) dw_s + \int_{t_n}^t \big( a y_1(s)
- a_2 \sqrt{y_1(s)} \big)ds, \\
 y_2(t)  =  y_2(t_n) & + & \int_{t_n}^t \Big(
  a y_1(s) - a_2
\sqrt{y_1(s)}+  a_1+y_2(s)\big(a_3-a+\sigma^2
y_2^{\rho-1}(\hat{s})-a_4y_2^{\frac{r-1}{2}}(\hat{s}) \big) ds \\
&  + & \int_{t_n}^t 2 \sigma y_2(s)
y_2^{\frac{\rho-1}{2}}(\hat{s}) dw_s + \int_{t}^{t_{n+1}}  a
y_1(s) - a_2 \sqrt{y_1(s)} ds
\end{eqnarray*}

\begin{proposition}
If $r+1 > 2 \rho$  then   we have the following moment bounds, for
$\Delta < 1$ if $a = \ln \frac{4}{3}$,
\begin{eqnarray*}
\mathbb{E}(|y_1(t)|^p+|y_2(t)|^p) < \infty
\end{eqnarray*}
\end{proposition}

\proof Let $\tau_R = \inf \{ t \in [0,T] : |y_2(t)| > R \}$. Note
that $y_1(t \wedge \tau_R)$ is also uniformly bounded for all
$\omega \in \Omega$.

 We can write
\begin{eqnarray*}
y_2(t \wedge \tau_R) = y_2(0) & + & \int_0^{t \wedge \tau_R} (a_1
- a_2 \sqrt{y_1(s )} + ay_1(s)+y_2(s)\big(a_3-a+\sigma^2
y_2^{\rho-1}(\hat{s})-a_4y_2^{\frac{r-1}{2}}(\hat{s}) \big) ds \\
& + & \int_{0}^{t\wedge \tau_R} 2 \sigma y_2(s)
y_2^{\frac{\rho-1}{2}}(\hat{s}) dw_s + \int_{t\wedge
\tau_R}^{t_{n+1}\wedge \tau_R}  ay_1(s) - a_2 \sqrt{y_1(s)} ds
\end{eqnarray*}
But
\begin{eqnarray*}
y_2(t\wedge \tau_R) \leq y_2(0)  & + & \int_0^{t\wedge \tau_R}
(a_1 + ay_1(s)+y_2(s)\big(a_3+\sigma^2
y_2^{\rho-1}(\hat{s})-a_4y_2^{\frac{r-1}{2}}(\hat{s}) \big) ds \\
& + & \int_{0}^{t\wedge \tau_R} 2 \sigma y_2(s)
y_2^{\frac{\rho-1}{2}}(\hat{s}) dw_s + \int_{t_n\wedge
\tau_R}^{t_{n+1}\wedge \tau_R} a y_1(s)  ds
\end{eqnarray*}
Let us estimate
\begin{eqnarray*}
\int_{t_n\wedge \tau_R}^{t_{n+1}\wedge \tau_R} a y_1(s)  ds & = &
a\int_{t_n\wedge \tau_R}^{t_{n+1}\wedge \tau_R}\Big(\frac{a_2}{a}
+ (\sqrt{y_2(t_n)} - \frac{a_2}{a}) e^{a\frac{(t-t_n)}{2}} \Big)^2
ds \\ & = &  a \int_{t_n\wedge \tau_R}^{t_{n+1}\wedge \tau_R}
\Big(\frac{a^2_2}{a^2} + (\sqrt{y_2(t_n)} - \frac{a_2}{a})^2
e^{a(t-t_n)} +2\frac{a_2}{a}(\sqrt{y_2(t_n)} -
\frac{a_2}{a})e^{a\frac{(t-t_n)}{2}} \Big) \\ & \leq & C
+\int_{t_n\wedge \tau_R}^{t_{n+1}\wedge \tau_R}a(\sqrt{y_2(t_n)} -
\frac{a_2}{a})^2 e^{a(t-t_n)} + 2a\frac{a^2_2}{a^2} + \frac{1}{2}a
(\sqrt{y_2(t_n)} - \frac{a_2}{a})^2 e^{a(t-t_n)}ds \\ & = & C +
\int_{t_n\wedge \tau_R}^{t_{n+1}\wedge
\tau_R}\frac{3}{2}a(\sqrt{y_2(t_n)} - \frac{a_2}{a})^2
e^{a(t-t_n)} \\ & \leq & C + \frac{3a_2^2\Delta}{2a} e^{a\Delta} +
\int_{t_n\wedge \tau_R}^{t_{n+1}\wedge \tau_R}\frac{3}{2}a
y_2(t_n)e^{a(t-t_n)}ds \\ & \leq & C + \frac{3}{2} y_2(t_n\wedge
\tau_R)(e^{a \Delta}-1)
\end{eqnarray*}

Therefore
\begin{eqnarray*}
y_2(t\wedge \tau_R) \leq  C & + & \frac{3}{2} y_2(t_n\wedge
\tau_R)(e^{a \Delta}-1)  + \int_0^{t\wedge \tau_R} (a_1  +
ay_1(s)+y_2(s)\big(a_3+\sigma^2
y_2^{\rho-1}(\hat{s})-a_4y_2^{\frac{r-1}{2}}(\hat{s}) \big) ds \\
 & + &  \int_{0}^{t\wedge
\tau_R} 2 \sigma y_2(s) y_2^{\frac{\rho-1}{2}}(\hat{s}) dw_s
\end{eqnarray*}

Denoting by $v(s\wedge \tau_R)$ the following Ito process,
\begin{eqnarray*}
v(s\wedge \tau_R) =  C & + & \frac{3}{2} y_2(t_n\wedge
\tau_R)(e^{a \Delta}-1) + \int_0^{t\wedge \tau_R} (a_1  +
ay_1(s)+y_2(s)\big(a_3+\sigma^2
y_2^{\rho-1}(\hat{s})-a_4y_2^{\frac{r-1}{2}}(\hat{s}) \big) ds \\
 & + &  \int_{0}^{t\wedge
\tau_R} 2 \sigma y_2(s) y_2^{\frac{\rho-1}{2}}(\hat{s}) dw_s
\end{eqnarray*}
we see that $y_2(s\wedge \tau_R) \leq v(s\wedge \tau_R)$. We shall
prove that $v(s)$ has bounded $pth$ moments and therefore $y_2(t)$
also.

 Applying
Ito's formula on $v^p(t \wedge \tau_R)$ we get that
\begin{eqnarray*}
v^p(t\wedge \tau_R) = & & (C+\frac{3}{2} y_2(t_n\wedge
\tau_R)(e^{a \Delta}-1) )^p \\ & + & \int_0^{t\wedge \tau_R}
pv^{p-1} \big(a_1 + ay_1(s)+y_2(s)\big(a_3+\sigma^2
y_2^{\rho-1}(\hat{s})-a_4y_2^{\frac{r-1}{2}}(\hat{s})\big)\big)ds
\\ & + &\int_0^{t\wedge \tau_R}\frac{p(p-1)}{2} 4 \sigma^2 y_2^{2}(s)
v^{p-2}(s)y_2^{\rho-1}(\hat{s}) ds \\ & + & \int_{0}^{t\wedge
\tau_R}p 2 \sigma y_2(s)v^{p-1}(s) y_2^{\frac{\rho-1}{2}}(\hat{s})
dw_s
\end{eqnarray*}
Taking expectations we arrive at
\begin{eqnarray*}
\mathbb{E}v^p(t \wedge \tau_R) \leq & &  \mathbb{E} (C+\frac{3}{2}
y_2(t_n\wedge \tau_R)(e^{a \Delta}-1) )^p \\ & & + \int_0^{t
\wedge \tau_R} ( C + C \mathbb{E}v^p(s) + pa\mathbb{E}y_1(s)v^{p-1}(s)ds \\
& & + \int_0^{t \wedge \tau_R}p\mathbb{E}v^{p}(s)
\big(a_3+\sigma^2
y_2^{\rho-1}(\hat{s})-a_4y_2^{\frac{r-1}{2}}(\hat{s}) +
\frac{p(p-1)}{2} 4 \sigma^2  y_2^{\rho-1}(\hat{s})\big)ds
\end{eqnarray*}
We have assumed that $r+1 > 2 \rho$ so there exists some constant
$C$ independent of $\omega,\Delta$ such that
\begin{eqnarray*}
a_3+\sigma^2 y_2^{\rho-1}(\hat{s})-a_4y_2^{\frac{r-1}{2}}(\hat{s})
+ \frac{p(p-1)}{2} 4 \sigma^2  y_2^{\rho-1}(\hat{s}) \leq C
\end{eqnarray*}

We will estimate now the term, using the inequality $(a+b)^p \leq
2^{p-1}a^p + 2^{p-1}b^b$,
\begin{eqnarray*}
\mathbb{E} (C+\frac{3}{2} y_2(t_n\wedge \tau_R)(e^{a \Delta}-1)
)^p \leq C + \frac{3^p}{2}  (e^{a \Delta}-1)^p
\mathbb{E}y_2^p(t_n\wedge \tau_R)
\end{eqnarray*}
Furthermore the following term
\begin{eqnarray*}
\mathbb{E}y_1(s)v^{p-1}(s) = \mathbb{E}\Big(\frac{a_2}{a} +
(\sqrt{y_2(\hat{s})} - \frac{a_2}{a}) e^{a\frac{(t-t_n)}{2}}
\Big)^2 v^{p-1}(s) \leq C + C\sup_{ 0 \leq l \leq
s}\mathbb{E}v^p(l)
\end{eqnarray*}
Collecting all the above results, we obtain
\begin{eqnarray*}
\mathbb{E}v^p(t \wedge \tau_R) \leq C + \frac{3^p}{2}  (e^{a
\Delta}-1)^p\mathbb{E}v^p(t_n\wedge \tau_R) + C\int_{0}^{t \wedge
\tau_R} \mathbb{E}v^p(s) + \sup_{ 0 \leq l \leq s}\mathbb{E}v^p(l)
ds
\end{eqnarray*}

 Next setting
\begin{eqnarray*}
u(s) = \sup_{ 0 \leq l \leq s} \mathbb{E} v^p(l)
\end{eqnarray*}
we can write
\begin{eqnarray*}
(1- \frac{3^p}{2}  (e^{a \Delta}-1)^p)u(t) \leq C + \int_{0}^{t
\wedge \tau_R} C u(s) ds
\end{eqnarray*}
Now it is the time to choose accordingly the free parameter $a$ so
as
\begin{eqnarray*}
 \frac{3^p}{2}  (e^{a \Delta}-1)^p < 1
 \end{eqnarray*}
 Choosing $a = \ln \frac{4}{3}$ we easily see that the above
 inequality holds for every $\Delta < 1$. It is clear that if we
 use smaller $\Delta$ then we can choose bigger $a$ so that the
 corresponding constants will be smaller.

Using Gronwall's inequality and then Fatou's lemma we get the
result.
 \qed

Unfortunately, we can not use our Theorem 1 directly to get the
desired result, therefore we will argue differently.

\begin{theorem} If $r+1 > 2 \rho$ then
\begin{eqnarray*}
\mathbb{E}|y(t) - y_2(t)|^2 \leq C_R \Delta + \frac{C}{R}
\end{eqnarray*}
for any $R > 0$. Therefore, for every $\vep > 0$ we fix  $R$ such
that $\frac{C}{R} < \vep$ and then for small enough $\Delta$ we
have that
\begin{eqnarray*}
\mathbb{E}|y(t) - y_2(t)|^2 \leq \vep
\end{eqnarray*}
\end{theorem}

\proof The approximate solution $y_2(t)$ is as follows,
\begin{eqnarray*}
 y_2(t)  =  x_0^2 & + & \int_0^t ( a_1 - a_2 \sqrt{y_1(s)} +
(a_3-a)y_2(s) + a y_1(s) - a_4 y_2(s)y_2^{\frac{r-1}{2}}(\hat{s})
+ \sigma^2 y_2(s)y_2^{\rho-1}(\hat{s}) ) ds \\ &  + & \int_0^t 2
\sigma y_2(s)y_2^{\frac{\rho-1}{2}}(\hat{s}) dw_s +
\int_t^{t_{n+1}} a y_1(s) - a_2 \sqrt{y_1(s)} ds
\end{eqnarray*}
Set
\begin{eqnarray*}
v(t)  =  x_0^2 & + & \int_0^t ( a_1 - a_2 \sqrt{y_1(s)} +
(a_3-a)y_2(s) + a y_1(s) - a_4 y_2(s)y_2^{\frac{r-1}{2}}(\hat{s})
+ \sigma^2 y_2(s)y_2^{\rho-1}(\hat{s}) ) ds \\ &  + & \int_0^t 2
\sigma y_2(s)y_2^{\frac{\rho-1}{2}}(\hat{s}) dw_s
\end{eqnarray*}
Then it is clear that $\mathbb{E}|y_2(t) - v(t)|^p \leq C
\Delta^p$ and therefore $v(t)$ has bounded moments as well.

The difference $y(t)-v(t)$ is as follows
\begin{eqnarray*}
y(t)-v(t) & = & \int_0^t \Big( a_2(\sqrt{y_1(s)} - \sqrt{y(s)})
+(a_3-a)(y(s)-y_2(s)) + a (y(s)-y_1(s))
  \\ & & +a_4(y_2(s)y_2^{\frac{r-1}{2}}(\hat{s}) - y^{\frac{r+1}{2}}(s)  +
\sigma^2 (y^{\rho}(s) - y_2(s)y_2^{\rho-1}(\hat{s})) \Big)ds \\ &
& + \int_0^t 2 \sigma (y^{\frac{\rho+1}{2}}(s) -
y_2(s)y^{\frac{\rho-1}{2}}_2(\hat{s})) dw_s
\end{eqnarray*}

Applying Ito's formula on $(y(t)-v(t))^2$ and setting
\begin{eqnarray*}
g(t)  & = & a_2(\sqrt{y_1(s)} - \sqrt{y(s)}) +(a_3-a)(y(s)-y_2(s))
+ a (y(s)-y_1(s))
   \\ & & +a_4(y_2(s)y_2^{\frac{r-1}{2}}(\hat{s}) - y^{\frac{r+1}{2}}(s))  +
\sigma^2 (y^{\rho}(s) - y_2(s)y_2^{\rho-1}(\hat{s}))
\end{eqnarray*}
 we obtain
\begin{eqnarray*}
\mathbb{E}(y(t)-v(t))^2 = \mathbb{E}\int_0^t 2 (y(t)-y_2(t)) g(s)
+ 2 (y_2(t)-v(t)) g(s)  + 4 \sigma^2 (y^{\frac{\rho+1}{2}}(s) -
y_2(s)y^{\frac{\rho-1}{2}}_2(\hat{s}))^2 ds
\end{eqnarray*}

Then, it is clear that
\begin{eqnarray*}
2 \mathbb{E} (y_2(t)-v(t)) g(s) \leq C \Delta
\end{eqnarray*}
using  the moment bounds of $y,y_2$.

Next we will estimate the term
\begin{eqnarray*}
2 a_2\mathbb{E} (y(t)-y_2(t)) (\sqrt{y_1(s)} - \sqrt{y(s)}) & = &
2 a_2\mathbb{E} (y(t)-y_2(t)) (\sqrt{y_2(s)} - \sqrt{y(s)}) \\ & &
+ 2 a_2\mathbb{E} (y(t)-y_2(t)) (\sqrt{y_1(s)} - \sqrt{y_2(s)}) \\
& \leq & 2 a_2\mathbb{E} (y(t)-y_2(t)) (\sqrt{y_1(s)} -
\sqrt{y_2(s)})
\end{eqnarray*}
where we have used the fact that $h(x) = -\sqrt{x}$ is a
decreasing function. Using Holder's inequality and the fact that
$y,y_1,y_2$ have bounded moments, we obtain
\begin{eqnarray*}
\mathbb{E} (y(t)-y_2(t)) (\sqrt{y_1(s)} - \sqrt{y_2(s)}) \leq
\sqrt{\mathbb{E} |y(t)-y_2(t)|^2} \sqrt{\mathbb{E}
|y_1(s)-y_2(s)|} \leq C \sqrt{\mathbb{E} |y_1(s)-y_2(s)|}
\end{eqnarray*}
Therefore it remains to estimate the term
\begin{eqnarray*}
\mathbb{E}|y_1(s) - y_2(s)| & \leq &
\mathbb{E}|y_1(\tilde{s})-y_2(s)| + \mathbb{E}|y_1(s) -
y_1(\tilde{s})| \\ & \leq & \mathbb{E}|y_1(\tilde{s})-y_2(s)| +
\mathbb{E}|y_1(\tilde{s})-y_1(\hat{s})| +
\mathbb{E}|y_1(\hat{s})-y_1(s)| \\ & \leq & C \sqrt{\Delta}
\end{eqnarray*}

Now we will estimate the term
\begin{eqnarray*}
& & 2 \mathbb{E} (y(t)-y_2(t))\Big((a_3-a)(y(s)-y_2(s)) + a
(y(s)-y_1(s))\Big) \\ &=& 2 a_3\mathbb{E} (y(t)-y_2(t))^2  +2 a
\mathbb{E} (y(t)-y_2(t)) (y_1(t)-y_2(t)) \\ & \leq & C\mathbb{E}
(y(t)-y_2(t))^2 + C \sqrt{\Delta}
\end{eqnarray*}

The term
\begin{eqnarray*}
& & 2 \mathbb{E} (y(t)-y_2(t))
(a_4(y_2(s)y_2^{\frac{r-1}{2}}(\hat{s}) - y^{\frac{r+1}{2}}(s)) \\
& = &  2 a_4\mathbb{E} (y(t)-y_2(t)) (y^{\frac{r+1}{2}}_2(s) -
y^{\frac{r+1}{2}}(s)) +2
a_4\mathbb{E}(y(t)-y_2(t))y_2(s)(y_2^{\frac{r-1}{2}}(\hat{s})-y_2^{\frac{r-1}{2}}(s))
\\ & \leq & C \Delta + 2
a_4\mathbb{E}(y(t)-y_2(t))y_2(s)(y_2^{\frac{r-1}{2}}(\hat{s})-y_2^{\frac{r-1}{2}}(s))
\end{eqnarray*}
using the fact that $h(x) = - x^{\frac{r+1}{2}}$ is a decreasing
function. But
\begin{eqnarray*}
& &
2a_4\mathbb{E}(y(t)-y_2(t))y_2(s)(y_2^{\frac{r-1}{2}}(\hat{s})-y_2^{\frac{r-1}{2}}(s))
 \\ & = & 2a_4\mathbb{E}(y(t)-y_2(t))\big( y_2^{\frac{r+1}{2}}(\hat{s}) -
 y_2^{\frac{r+1}{2}}(s)\big) +
 2a_4\mathbb{E}(y(t)-y_2(t))(y_2(s)-y_2(\hat{s}))y_2^{\frac{r-1}{2}}(\hat{s})
 \end{eqnarray*}
Since $\mathbb{E}y^{\frac{r-1}{2}}_2(t) < \infty$ we can use the
mean value theorem arriving to the following estimate,
\begin{eqnarray*}
\mathbb{E} (y(t)-y_2(t)) (a_4(y_2(s)y_2^{\frac{r-1}{2}}(\hat{s}) -
y^{\frac{r+1}{2}}(s)) \leq C \sqrt{\Delta}
\end{eqnarray*}

The term
\begin{eqnarray*}
2 \mathbb{E}(y(t)-y_2(t))(y^{\rho}(t) - y_2^{\rho}(t)) = 2(\rho-1)
\mathbb{E} (y(t)-y_2(t))^2 h^{\rho-1}(t)
\end{eqnarray*}
using the mean value theorem, for some $h(t)$ located between
$y(t),y_2(t)$. Using the moment bounds of $y,y_2$ we arrive at,
for any $R>0$,
\begin{eqnarray*}
2(\rho-1) \mathbb{E} (y(t)-y_2(t))^2 h^{\rho-1}(t) & = & 2(\rho-1)
\mathbb{E} (y(t)-y_2(t))^2 h^{\rho-1}(t) \mathbb{I}_{ \{
h^{\rho-1}(t) > R \}} \\ & & + 2(\rho-1) \mathbb{E}
(y(t)-y_2(t))^2 h^{\rho-1}(t) \mathbb{I}_{ \{ h^{\rho-1}(t) \leq R
\}} \\ & \leq & C_R \mathbb{E} |y(t)-y_2(t)|^2  +2(\rho-1)
\mathbb{E} (y(t)-y_2(t))^2 h^{\rho-1}(t) \mathbb{I}_{ \{
h^{\rho-1}(t) > R \}}
\end{eqnarray*}
Using Holder's inequality we deduce that
\begin{eqnarray*}
\mathbb{E} (y(t)-y_2(t))^2 h^{\rho-1}(t) \mathbb{I}_{ \{
h^{\rho-1}(t) > R \}} \leq \sqrt{\mathbb{E} (y(t)-y_2(t))^4
h^{2\rho-2}(t)} \sqrt{\mathbb{E} \mathbb{I}_{ \{ h^{\rho-1}(t) > R
\}}}
\end{eqnarray*}
Using the moment bounds of $y,y_2,h$ and Markov's inequality we
arrive at the following estimate,
\begin{eqnarray*}
2 \mathbb{E}(y(t)-y_2(t))(y^{\rho}(t) - y_2^{\rho}(t)) \leq C_R
\mathbb{E} |y(t)-y_2(t)|^2 + \frac{C}{R}
\end{eqnarray*}

The same holds for the term $4 \sigma^2
\mathbb{E}(y^{\frac{\rho+1}{2}}(s) -
y_2(s)y^{\frac{\rho-1}{2}}_2(\hat{s}))^2$ therefore we conclude
that
\begin{eqnarray*}
\mathbb{E}(y(t)-v(t))^2 \leq \frac{C}{R} + C \sqrt{\Delta} +
C_R\int_0^t \mathbb{E}(y(s)-v(s))^2 ds
\end{eqnarray*}
and using Gronwall's inequality we conclude the desired result.

Now it is easy to see that
\begin{eqnarray*}
\mathbb{E}|y(t)-y_2(t)|^2 \leq C \mathbb{E}|y(t)-v(t)|^2 + C
\mathbb{E}|y_2(t)-v(t)|^2 \to 0 \mbox{ as } \Delta \to 0
\end{eqnarray*}
 \qed

 \begin{theorem} If $r+1> 2 \rho$ we have that
 \begin{eqnarray*}
\lim_{\Delta \to 0} \mathbb{E} | \sqrt{y_2(t)} - \sqrt{y(t)} |^2 =
0
\end{eqnarray*}
\end{theorem}

\proof Using Holder's inequality and the nonnegativity of the
approximation we have
\begin{eqnarray*}
\mathbb{E} |\sqrt{y_2(t)} - \sqrt{y(t)} |^2 & = & \mathbb{E}
\frac{|y_2(t) - y(t)|^2}{(\sqrt{y(t)} + \sqrt{y_2(t)})^2} \\ &
\leq & \sqrt{\mathbb{E}|y_2(t) - y(t)|^2}
\sqrt[4]{\mathbb{E}|y(t)+y_2(t)|^4\mathbb{E}\frac{1}{y^2(t)}}
\end{eqnarray*}
Recalling from \cite{Mao} that the inverse moments of the true
solution are bounded we conclude our proof. \qed

In practice, we will approximate equation (8) with a slightly
different approximation which is the following
\begin{eqnarray*}
y_2(t) = y_1(t_{n+1}) + a_1 \Delta+ \int_{t_n}^t
y_2(s)\big(\sigma^2
y_2^{\rho-1}(\hat{s})-a_4y_2^{\frac{r-1}{2}}(\hat{s}) \big) ds +
\int_{t_n}^t 2 \sigma y_2(s) y_2^{\frac{\rho-1}{2}}(\hat{s}) dw_s
\end{eqnarray*}

Our numerical scheme for the transformed Ait-Sahalia model is as
follows
\begin{eqnarray*}
y_0 & = & x_0^2, \\
y_{n+1} &  = & \Big(a_1\Delta + \left(\frac{a_2}{a} +
(\sqrt{y_n}-\frac{a_2}{a}) e^{\frac{a\Delta}{2}}\right)^2 \Big)
e^{-\Delta(\sigma^2y_n^{\rho-1} + a_4 y_n^{\frac{r-1}{2}})  + 2
\sigma y_n^{\frac{\rho-1}{2}} \Delta W}
\end{eqnarray*}

The approximation of the Ait-Sahalia model will be $x_n =
\sqrt{y_n}$.

\section{Summary and comments}
In this paper we extend the semi discrete method by combine it
with the split step method. We can, in general, split our
stochastic differential equation in $m$ equations (by splitting
the drift term) and then in each equation we apply any
approximation method. The aim of this new technique is that the
resulting numerical scheme will be boundary preserving.

Consider for example the following Ait-Sahalia type model,
\begin{eqnarray*}
x_t = x_0 + \int_0^t \frac{1}{2}\Big(\frac{a_1}{x_s} - a_2 -b_1
\sqrt{x_s} - b_2 x^{\frac{3}{2}}_s - b_3 \ln (1+x_s^2) + a_3 x_s -
a_4 x^{r}_s \Big) ds + \int_0^t \sigma x_s^{\rho} dw_s.
\end{eqnarray*}
Using the transformation $y_t = x_t^2$ we get the following
stochastic differential equation for $y_t$,
\begin{eqnarray*}
y_t = y_0 & + &  \int_0^t \Big(a_1-a_2 \sqrt{y_s} - b_1
y_s^{\frac{3}{4}} - b_2 y^{\frac{5}{4}}_s - b_3 \sqrt{y_s} \ln
(1+y_s) + a_3 y_s - a_4 y_s^{\frac{r+1}{2}}  +  \sigma^2
y_s^{\rho} \Big) ds \\ & + & \int_0^t 2 \sigma
y_s^{\frac{\rho+1}{2}} dw_s
\end{eqnarray*}
Then we propose the following slit step combined with the semi
discrete technique numerical scheme, for $t \in ( t_n, t_{n+1}]$
\begin{eqnarray*}
y_5(0) & = & x_0,  \\
y_1(t) & =& y_5(t_n) -b_3 \int_{t_n}^t \sqrt{y_1(s)} \ln
(1+y_5(t_n))
ds, \\
y_2(t) & = & y_1(t_{n+1}) -b_2 \int_{t_n}^t
y_2(s)y_5^{\frac{1}{4}}(t_n) ds,
\\
y_3(t) & = & y_2(t_{n+1}) - b_1 \int_{t_n}^t y_3^{\frac{3}{4}}(s)
ds,
\\
y_4(t) & = & y_3(t_{n+1}) + \int_{t_n}^t a_3 y_4(s) -
a_2\sqrt{y_4(s)}
ds, \\
y_5(t) & = & y_4(t_{n+1}) + a_1\Delta+ \int_{t_n}^t y_5(s)
\big(-a_4 y_5^{\frac{r-1}{2}}(t_n) + \sigma^2 y_5^{\rho-1}(t_n)
\big) ds + \int_{t_n}^t 2 \sigma y_5^{\frac{\rho-1}{2}}(t_n)y_5(s)
dw_s
\end{eqnarray*}
The solutions are
\begin{eqnarray*}
y_1(t) & = & \Big(\sqrt{y_5(t_n)}-\frac{b_3
\ln(1+y_5(t_n))}{2}(t-t_n) \Big)^2, \\
y_2(t) & = & y_1(t_{n+1}) e^{-b_2y_5(t_n) (t-t_n)}, \\
y_3(t) & = & \Big(\sqrt[4]{y_2(t_{n+1})} - b_1\frac{t-t_n}{4}
\Big)^4,\\
y_4(t) & = & \Big(\frac{a_2}{a_3} + (\sqrt{y_3(t_n)} -
\frac{a_2}{a_3}) e^{a_3\frac{(t-t_n)}{2}} \Big)^2, \\
y_5(t) & = & (y_4(t_{n+1}) + a_1\Delta)
e^{-\Delta(\sigma^2y_5^{\rho-1}(t_n) + a_4
y_5^{\frac{r-1}{2}}(t_n)) + 2 \sigma y_5^{\frac{\rho-1}{2}}(t_n)
\Delta W}
\end{eqnarray*}

This use of the splitting-semi discrete technique produces a
positivity preserving and explicit numerical scheme.

Another, obviously generalization, is that we can semi-discretize
in the time variable also. For example we can assume the following
assumption,
\\[0.2cm]
{\bf Assumption C} Suppose that the functions $f_i$ for
$i=1,...,m$ satisfy the following locally Lipschitz condition,
\begin{eqnarray*}
 |f_i(t,t,x,x,x) - f_i(t,t_1,x_1,x_2,x_3)|  & \leq &  C_R (|t-t_1|^{b_1}+|x-x_1| +
|x-x_2|+|x-x_3| ), \quad i=2,...,m \\
|f_{m+1}(t,t,x,x) - f_{m+1}(t,t_1,x_1,x_2)| & \leq & C_R
(|t-t_1|^{b_2}+|x-x_1| + |x-x_2| + |x-x_1|^a)
\end{eqnarray*}
for $|x_1|, |x_2|, |x_3|, |x| \leq R$, for $t,t_1 \in [0,T]$, for
some $a \in [\frac{1}{2},1)$ and for some $b_1,b_2 > 0$. With this
kind of generalization we can handle problems which has time
depending parameters (see \cite{Halidias1}).

There is also the possibility to split the diffusion term. For
example, if $a_1+...+a_m = a$ and $b_1+...+b_l = b$ then for $ t
\in (t_n,t_{n+1}]$,
\begin{eqnarray*}
y_m(0) & = & x_0, \nonumber \\
y_1(t) & = & y_m(t_n) + \int_{t_n}^t a_1(s,y_1(s)) ds, \\
 & \vdots & \\
 y_i(t) & = & y_{i-1}(t_{n+1}) + \int_{t_n}^t a_i(s,y_i(s))ds +
 \int_{t_n}^t b_1(s,y_i(s))dw_s, \\
& \vdots & \\
 y_j(t) & = & y_{j-1}(t_{n+1}) + \int_{t_n}^t a_j(s,y_j(s))ds, \\
& \vdots & \\
y_m(t) & =  & y_m(t_{n+1}) + \int_{t_n}^t a_m(s,y_2(s)) ds +
\int_{t_n}^t b_{l}(s,y_m(s))dw_s
\end{eqnarray*}
where $i,j \in \{1,...,m \}$. The general idea is that if we want
to construct a numerical scheme with values in a domain $D$ then
we can split accordingly such as $y_1 \in D_1$, whenever the
initial value (for $y_1$) is in $D$, $y_2 \in D_2$ whenever the
initial value is in $D_1$ and finally $y_m \in D$ whenever the
initial value for $y_m$ is in $D_{m-1}$.

To approximate any of the above equations which do not have a
diffusion term we can use any suitable numerical scheme and any
semi discrete approximation. The same holds also for the first
stochastic differential equation (in our setting is the
$i$-equation) but to approximate any other stochastic differential
equation (i.e. with a diffusion term) we should fully discretize
the corresponding sde, i.e. we can not use a semi discrete method.
If we want to produce a boundary preserving numerical scheme then
these  sdes (that we should fully discretize) can be approximated
by balanced Milstein methods (see \cite{Kahl}, \cite{Milstein},
\cite{Milstein2} and \cite{Dan}).

An interesting question (that we do not answer in this paper) is
the rate of convergence of the explicit numerical schemes that the
semi discrete method produces. In  \cite{Hunt1}, \cite{Hunt2},
\cite{Jentzen}  the authors study the convergence rates of various
numerical schemes and it seems that these techniques will be
useful also for the semi discrete schemes.

\end{document}